\documentclass[12pt,reqno]{amsart}

\usepackage{amsmath, amsfonts, amsthm, amssymb, amscd, dsfont, tipa,mathtools, listings, hyperref}

\newtheorem{thm}{Theorem}[section]

\newtheorem{cor}[thm]{Corollary}
\newtheorem{lem}[thm]{Lemma}
\newtheorem{prop}[thm]{Proposition}

\newtheorem{rek}[thm]{Remark}

\newcommand{\PSU}{\text{PSU}}
\newcommand{\Aut}{\text{Aut}}
\newcommand{\Ppd}[2]{\mathcal P_{#1}^{#2}}
\newcommand{\ord}{\text{ord}}
\newcommand{\SU}{\text{SU}}
\newcommand{\PGU}{\text{PGU}}

\title{The residual unitary case for almost elusive groups}
\author{Alan Xuelun Hou} \address[Alan Hou]{Brandeis University, Waltham, MA 02453}\email{hou@brandeis.edu}
\author{Tudor Popescu} \address[Tudor Popescu]{Brandeis University, Waltham, MA 02453}\email{tudorpopescu@brandeis.edu}
\date{}

\begin{document}

\begin{abstract}
    In this short note, we remove the residual unitary family in Hall's classification of almost elusive classical groups and the corresponding conditional family in the quasiprimitive classification. Using theorems of Zsigmondy and Schinzel together with some elementary number theory arguments, we prove that if $n \ge 5$ is odd, $f \ge 2$ is even, and $r$ is the unique primitive prime divisor of $(2^f)^{2n} - 1,$ then $4nf \mid r - 1.$ For $n \ge 5$ prime and any integer $f \ge 1,$ we deduce that a unique primitive prime divisor $r$ of $(2^f)^{2n} - 1$ satisfies $r \ge 4nf + 1,$ except when $(n, f, r) = (5, 1, 11).$ We then obtain two conjugacy classes of derangements of the same prime order in a broader family of unitary subspace actions, and a quadratic lower bound for the number of such classes in the corresponding product actions.
\end{abstract} 
\maketitle

\section{Introduction} 
A derangement in a finite transitive permutation group is an element that fixes no point. In \cite{BurnessHall2022}, Burness and Hall define an \emph{almost elusive group} as a permutation group that has exactly one conjugacy class of derangements of prime order. Hall's classification of primitive almost elusive almost simple classical groups leaves one residual family: the socle is $\PSU_n(q),$ where $q$ is even and $n \ge 5$ is a prime divisor of $q + 1,$ and the point stabilizer stabilizes a nondegenerate $1$-space \cite[Theorem 1(i)]{Hall2022}. This is the only potentially infinite family in the range of that theorem, and Hall anticipates that it is empty \cite[Remarks 1 and 5.34]{Hall2022}. We prove this is the case and obtain a stronger statement for unitary subspace actions.

For integers $a > 1$ and $m \ge 2,$ let
$$\Ppd{a}{m} = \{r \text{ prime} : \ord_r(a)=m\}$$
be the set of primitive prime divisors of $a^m - 1.$

\begin{thm}\label{thm:main}
Let $f \ge 1$ be an integer, set $q = 2^f,$ and let $n \ge 5$ be odd. Set $T = \PSU_n(q)$ and suppose that either $f$ is even or $n \ge 7$ is prime. Then there is a prime $r \in \Ppd{q}{2n}$ such that the elements of order $r$ in $T$ form at least two $\Aut(T)$-classes.
\newline
For every $T \le G \le \Aut(T)$, all these elements are derangements in the natural action of $G$ on any of its orbits of proper nonzero subspaces of the defining Hermitian space over $\mathbb{F}_{q^2}$. In particular, each such action has at least two $G$-classes of derangements of the same prime order.
\end{thm}

\begin{cor}\label{cor:classification}
Let $q=2^f$, where $f \ge 1$ is an integer, and let $n \ge 5$ be a prime divisor of $q+1$. Set $T=\PSU_n(q)$ and let $T \le G \le \Aut(T)$ act naturally on the nondegenerate $1$-spaces of the defining Hermitian space. Then $G$ has at least two conjugacy classes of derangements of the same prime order. Consequently, case (i) of \cite[Theorem 1]{Hall2022}, case U1 of \cite[Table A1]{Hall2024}, and the corresponding case (i) of \cite[Theorem 4]{Hall2024} do not occur.
\end{cor}

In \cite[Lemma 2.12]{Hall2022}, Hall leaves open the possibility
$$\Ppd{q}{2n} = \{2nf + 1\}$$
when $q = 2^f$ and $n \ge 5$ is a prime divisor of $q + 1.$ This is the necessary arithmetic condition in the residual case of \cite[Proposition 5.33]{Hall2022}.

In Proposition \ref{prop:divisibility}, we prove that if $n \ge 5$ is odd, $f \ge 2$ is even, and $\Ppd{2^f}{2n} = \{r\},$ then
$$4nf \mid r - 1.$$
Combining this with Hall's reduction, Corollary \ref{cor:lower-bound} gives $r \ge 4nf + 1$ whenever $n \ge 5$ is prime and $\Ppd{2^f}{2n} = \{r\}$, apart from $(n, f, r) = (5, 1, 11).$ Section 3 records the exact number of $\Aut(T)$-classes of irreducible elements of these prime orders and proves Theorem \ref{thm:main}. Section 4 applies this result to obtain a lower bound for the number of classes of prime-order derangements in product actions.

\section{Primitive prime divisors}
We use the following special cases of two classical theorems
on primitive prime divisors.

\begin{thm}[Zsigmondy {\cite{Zsigmondy1892}}]
Let $M \ge 2$ be an integer. If $M \neq 6$, then
$$\Ppd{2}{M} \ne \emptyset.$$
\end{thm}

\begin{thm}[Schinzel {\cite{Schinzel1962}}]
Let $M \ge 4$ be an integer with $M \equiv 4 \pmod{8}$.
If $M \notin \{4,12,20\}$, then
$$\bigl|\Ppd{2}{M}\bigr| \ge 2.$$
\end{thm}

We first prove the divisibility result for even $f$.

\begin{prop}\label{prop:divisibility}
Let $n \ge 5$ be odd, and let $f \ge 2$ be even. If
$$\Ppd{2^f}{2n} = \{r\},$$
then $4nf \mid r - 1.$
\end{prop}

\begin{proof}
Let $N = 2nf$. Since $N \ge 20$, Zsigmondy's theorem implies that $\Ppd{2}{N}\neq\emptyset,$ and let $s\in\Ppd{2}{N}.$ We have that
$$\ord_s(2^f) = \frac{\ord_s(2)}{\gcd(\ord_s(2), f)} = \frac{N}{f}=2n$$
As $\Ppd{2^f}{2n}=\{r\},$ we must have that $s = r.$ Therefore, $\ord_r(2) = N.$
\newline
If $v_2(f) = 1,$ then $N \equiv 4 \pmod{8}$. Every prime in $\Ppd{2}{N}$ belongs to $\Ppd{2^f}{2n},$ so Schinzel's theorem forces $N \in \{4,12,20\}.$ Since $N \ge 20,$ we get $(n, f) = (5, 2).$ In this case, $\Phi_{10}(4) = 205 = 5\cdot41$ and $\Ppd{4}{10} = \{41\},$ where $\Phi_{10}$ is the tenth cyclotomic polynomial, and the conclusion holds.
\newline
Otherwise, $4 \mid f$. As $\ord_r(2) = N$, we have $N \mid r - 1.$ Since $8 \mid N$, we get $r \equiv 1 \pmod{8}$. Therefore, $\left(\frac{2}{r}\right) = 1,$ and Euler's criterion yields
$$2nf = \ord_r(2) \mid \frac{r - 1}{2}$$
Thus $4nf \mid r - 1,$ as desired.
\end{proof}

\begin{cor}\label{cor:lower-bound}
Let $n \ge 5$ be prime, and let $f \ge 1$ be an integer. If $\Ppd{2^f}{2n}=\{r\}$, then
$$r \ge 4nf + 1$$
except when $(n, f, r) = (5, 1, 11).$
\end{cor}

\begin{proof}
By \cite[Lemma 2.12]{Hall2022}, either $r \ge 4nf + 1,$ or $r = 2nf + 1$ and one of the following holds: $(n, 2^f) = (5, 2)$ or $n \mid 2^f+1.$ The same lemma gives $f=2^a n^b$ for some integers $a,b \ge 0.$
\newline
The case $(n, 2^f) = (5, 2)$ gives the stated exception. Now suppose that $n \mid 2^f+1$. If $f$ is odd, then $f = n^b,$ and Fermat's little theorem gives
$$2^f + 1 = 2^{n^b} + 1 \equiv 3 \pmod{n}$$
contrary to $n \ge 5$ and $n \mid 2^f+1$. Therefore, $f$ is even. Proposition \ref{prop:divisibility} now gives $4nf \mid r - 1,$ which contradicts $r = 2nf + 1.$
\end{proof}

\begin{rek}
The lower bound is sharp, as $\Ppd{4}{10} = \{41\}$ and $41 = 4\cdot5\cdot2 + 1.$ The assumption that $f$ is even in Proposition \ref{prop:divisibility} cannot be omitted: indeed, $\Ppd{2}{14}=\{43\}$, but $28 \nmid 42$.
\end{rek}

\section{Prime order derangements in unitary groups}
We first record the following specialization of the class fusion discussion in \cite[Section 3.2.2]{Hall2022}.

\begin{lem}\label{lem:classes}
Let $n \ge 5$ be odd, $q = 2^f$ with $f \ge 1,$ and set $T = \PSU_n(q).$ For $r \in \Ppd{q}{2n},$ let $c_r$ denote the number of $\Aut(T)$-classes of elements of order $r$ in $T.$ Then
$$c_r = \frac{r - 1}{\ord_r(2)}$$
\end{lem}

\begin{proof}
Since $\ord_r(q)=2n,$ we have
$$\ord_r(q^2) = n, \qquad q^n \equiv -1 \pmod{r}, \qquad r \equiv 1 \pmod{2n}.$$
In particular, $r > n \ge |Z(\SU_n(q))|.$ Every element of order $r$ in $T$ therefore has a lift of order $r$ to $\SU_n(q).$ Such a lift is semisimple, and any nontrivial eigenvalue has degree $n$ over $\mathbb{F}_{q^2}.$ Therefore, the lift is irreducible.
\newline
The nontrivial $r$th roots of unity have $q^2$-Frobenius orbits of length $n.$ Since $q^n \equiv -1 \pmod{r},$ these orbits parametrize the irreducible conjugacy classes of elements of order $r$ in $\mathrm{GU}_n(q).$ Every such element lies in $\SU_n(q),$ since its determinant has order dividing both $r$ and $q + 1.$ Passing to the projective quotient introduces no further identifications: if $A^g = zB,$ where $A,B \in \SU_n(q)$ have order $r, g \in \mathrm{GU}_n(q),$ and $z \in Z(\SU_n(q)),$ then $z^r = 1,$ so $z = 1.$ Consequently, the elements of order $r$ in $T$ form $(r - 1)/n$ classes under $\PGU_n(q).$
\newline
The group $\Aut(T)$ is generated by $\PGU_n(q)$ and the field automorphism induced by squaring matrix entries, which squares eigenvalues. Since $q^2 = 2^{2f},$ the $\Aut(T)$-classes are therefore parametrized by the orbits of multiplication by $2$ on the nonzero exponents modulo $r.$ Every such orbit has length $\ord_r(2),$ so there are exactly $(r - 1)/\ord_r(2)$ classes.
\end{proof}

We are now ready to prove the main result.

\begin{proof}[Proof of Theorem \ref{thm:main}]
First suppose that $f$ is even, and set $N = 2nf.$ Every prime $r \in \Ppd{2}{N}$ belongs to $\Ppd{q}{2n},$ and Lemma \ref{lem:classes} gives
$$c_r=\frac{r - 1}{N}$$
If $4 \mid f$, choose any $r \in \Ppd{2}{N}.$ Such an $r$ exists by Zsigmondy's theorem. Since $8 \mid N$ and $N \mid r - 1,$ we have $r \equiv 1 \pmod{8}.$ As in the proof of Proposition \ref{prop:divisibility}, Euler's criterion gives $N \mid (r - 1)/2,$ so $c_r \ge 2.$
\newline
If $v_2(f) = 1,$ then $N \equiv 4 \pmod{8}$ and $N \ge 20.$ When $N = 20,$ take $r = 41,$ so $c_r = 2.$ Otherwise, Schinzel's theorem gives two distinct primes in $\Ppd{2}{N}.$ Their indices $(r - 1)/N$ are distinct positive integers, so at least one is at least $2.$ This proves the assertion when $f$ is even.
\newline
Assume now that $f$ is odd and $n \ge 7$ is prime. Let $f = tu, t = n^{v_n(f)}, \gcd(u,2n) = 1,$ and set $Q = \Ppd{2^t}{2n}.$ We claim that
$$Q = \Ppd{2}{2nt} \subseteq \Ppd{q}{2n}$$
Indeed, if $s \in Q$ and $h = \ord_s(2),$ then $h/\gcd(h,t) = 2n.$ Since $t$ is a power of the prime $n,$ this forces $h = 2nt.$ Conversely, if $s \in \Ppd{2}{2nt},$ then
$$\ord_s(2^t) = \frac{2nt}{\gcd(2nt,t)}
= 2n$$
so $s \in Q.$ Moreover, if $s \in Q,$ then
$$\ord_s(q) = \ord_s((2^t)^u) = \frac{2n}{\gcd(2n,u)} = 2n$$
which proves the inclusion.
\newline
Zsigmondy's theorem gives $Q \neq \emptyset.$ If $Q$ contains two distinct primes, their indices $(r - 1)/(2nt)$ are distinct positive integers, so one is at least $2.$ If $Q = \{r\},$ then Corollary \ref{cor:lower-bound}, applied with $t$ in place of $f,$ gives $r \ge 4nt + 1,$ since $n \ge 7$ excludes its exception. In either case, there is a prime $r \in Q$ for which
$$c_r = \frac{r - 1}{\ord_r(2)} = \frac{r - 1}{2nt} \ge 2$$
Finally, as shown in the proof of Lemma \ref{lem:classes}, every element of order $r$ in $T$ has an irreducible lift to $\SU_n(q).$ It therefore fixes no proper nonzero subspace and is a derangement in each of the actions in the statement. Distinct $\Aut(T)$-classes remain distinct under every subgroup $T \le G \le \Aut(T),$ which completes the proof.
\end{proof}

\begin{proof}[Proof of Corollary \ref{cor:classification}]
If $n \ge 7,$ Theorem \ref{thm:main} applies directly. If $n = 5,$ then $5 \mid 2^f + 1$ forces $f \equiv 2 \pmod{4},$ so the theorem again applies. This gives two $G$-classes of derangements of the same prime order on the nondegenerate $1$-spaces and excludes the residual primitive family.
\newline
Let $M$ be the stabilizer in $G$ of a nondegenerate $1$-space. If $H < M,$ then every derangement on $G/M$ remains a derangement on $G/H.$ Therefore, the same two classes exclude the corresponding quasiprimitive family.
\end{proof}

\begin{rek}
Corollary \ref{cor:classification} excludes case (i) of \cite[Theorem 1]{Hall2022}, as well as case U1 of \cite[Table A1]{Hall2024}; see \cite[Remark 6.2(b)]{Hall2024}. In particular, U1 is empty. The unitary alternatives in \cite[Corollary 5(ii)(a) and Corollary 6(i)]{Hall2024} are empty as well.
\end{rek}

\begin{rek}\label{rem:symplectic}
Let $n$ and $q = 2^f$ satisfy the hypotheses of Theorem \ref{thm:main}, and set $S=\mathrm{Sp}_{2n}(q).$ Choose $r \in \Ppd{q}{2n}$ as in that theorem. Then the elements of order $r$ in $S$ are irreducible and form exactly
$$\frac{r - 1}{\ord_r(2)} \ge 2$$
$\Aut(S)$-classes. Indeed, their $S$-classes are parametrized by the $q$-Frobenius orbits on the nontrivial $r$th roots of unity \cite[Proposition 2.9(2)]{BreuerGuralnickKantor2008}, and the field automorphisms fuse these classes by squaring eigenvalues. Since $\ord_r(q) = 2n,$ we have $r\nmid | \mathrm{O}_{2n}^{+}(q)|.$ Thus these elements are derangements in the action of $S$ on the plus-type quadratic forms with the prescribed symplectic polar form. Consequently, every group $S \le G \le \Aut(S)$ has at least two $G$-classes of derangements of the same prime order in this action.
\end{rek}

\section{Product actions}
For a finite transitive permutation group $W$ and a prime $r,$ let $K_r(W)$ denote the number of $W$-classes of derangements of order $r.$ We finish with a consequence for product actions.

\begin{prop}\label{prop:product}
Let $T$ be a finite nonabelian simple group acting transitively on a finite set $\Gamma,$ and suppose that this action extends to $\Aut(T).$ Suppose that $T$ contains $c \ge 1$ distinct $\Aut(T)$-classes of derangements on $\Gamma,$ all of the same prime order $r.$ Let $k \ge 1,$ and let
$$T^k \le W \le \Aut(T)\wr S_k = \Aut(T)^k \rtimes S_k,$$
where $S_k$ acts by permuting the $k$ factors of $\Aut(T)^k.$ Suppose that $W$ acts on $\Gamma^k$ with its product action. Then
$$K_r(W) \ge \binom{k+c}{c} - 1$$
In particular, if $T = \PSU_n(2^f)$ satisfies the hypotheses of Theorem \ref{thm:main} and $\Gamma$ is the set of nondegenerate $1$-spaces, then there is a prime $r$ such that
$$K_r(W) \ge \binom{k+2}{2} - 1 = \frac{k(k+3)}{2}$$
\end{prop}

\begin{proof}
Choose representatives $x_1, \ldots, x_c$ of the given $\Aut(T)$-classes. Consider the elements of $T^k$ whose coordinates belong to $\{1, x_1, \ldots, x_c\},$ excluding the identity. Every such element has order $r.$ It is also a derangement on $\Gamma^k,$ since at least one coordinate is a derangement on $\Gamma.$
\newline
Conjugation in $\Aut(T)\wr S_k$ permutes the coordinates and applies automorphisms of $T$ to them. It therefore preserves the multiplicity vector
$$(a_0, a_1, \ldots, a_c), \qquad a_0 + \ldots + a_c = k,$$
where $a_0$ counts identity coordinates and $a_i$ counts coordinates belonging to the class of $x_i.$ Distinct multiplicity vectors give elements that are not conjugate in the full wreath product, and hence are not conjugate in $W.$ There are $\binom{k+c}{c}$ such vectors of nonnegative integers. Removing the vector corresponding to the identity gives the required bound, and so the final assertion follows from Theorem \ref{thm:main} with $c = 2.$
\end{proof}

\begin{rek}
The bound is attained for a fixed prime order. For example take $T = \PSU_5(4)$ and $r = 41.$ Since $\ord_{41}(2) = 20,$ Lemma \ref{lem:classes} yields $c_{41} = 2.$ Let $W = \Aut(T)\wr S_k,$ with its product action on the nondegenerate $1$-spaces in each coordinate, and suppose $1 \le k < 41.$ Since $|\Aut(T):T|=20,$ the quotient $W/T^k$ has order $20^k k!,$ which is not divisible by $41.$ Therefore, we have that every element of order $41$ in $W$ lies in $T^k.$ Every such element is a derangement, and its coordinates belong to the identity class or one of the two $\Aut(T)$-classes of elements of order $41.$ Two such tuples are conjugate in the full wreath product if and only if they have the same multiplicity vector. Therefore,
$$K_{41}(W) = \binom{k+2}{2} - 1 = \frac{k(k+3)}{2}$$
This equality concerns the fixed prime order $41$. There are also derangements of order $5$ in $\PGU_5(4)\setminus T;$ see \cite[proof of Proposition 5.33]{Hall2022}. In particular, the total number of classes of derangements of prime order in $W$ is strictly larger than $K_{41}(W).$
\end{rek}

\end{document}